\documentclass[12pt,reqno]{amsart}
\usepackage{enumerate, latexsym, amsmath, amsfonts, amssymb, amsthm, color}
\def\pmod #1{\ ({\rm{mod}}\ #1)}

\def\N{\Bbb N}

\def\l{\left}
\def\r{\right}
\def\bg{\bigg}
\def\({\bg(}
\def\){\bg)}
\def\t{\text}
\def\f{\frac}

\def\ord{{\rm ord}}

\def\ls{\leqslant}

\def\bi{\binom}

\def\eq{\equiv}

\def\da{\delta}

\def\Proof{\noindent{\it Proof}}
\def\Ack{\medskip\noindent {\bf Acknowledgment}}
\theoremstyle{plain}
\newtheorem{theorem}{Theorem}

\newtheorem{lemma}{Lemma}

\theoremstyle{definition}

\theoremstyle{remark}
\newtheorem{remark}{Remark}

 \vspace{4mm}

\begin{document}

\hbox{Preprint}

\title
[{Congruences involving products of binomial coefficients}]
{New congruences involving products of two binomial coefficients}

\author
[Guo-Shuai Mao and Zhi-Wei Sun] {Guo-Shuai Mao and Zhi-Wei Sun}

\address {(Guo-Shuai Mao) Department of Mathematics, Nanjing
University, Nanjing 210093, People's Republic of China}
\email{mg1421007@smail.nju.edu.cn}

\address{(Zhi-Wei Sun) Department of Mathematics, Nanjing
University, Nanjing 210093, People's Republic of China}
\email{zwsun@nju.edu.cn}

\keywords{Central binomial coefficients, congruences.
\newline \indent The second author is the corresponding author. This research was supported by the Natural Science Foundation of China (grant 11571162).}

\begin{abstract} Let $p>3$ be a prime and let $a$ be a positive integer. We show that if $p\equiv1\pmod 4$ or $a>1$ then
$$\sum_{k=0}^{\lfloor\frac34p^a\rfloor}\frac{\binom{2k}k^2}{16^k}\equiv\l(\frac{-1}{p^a}\r)\pmod{p^3}$$
with $(-)$ the Jacobi symbol, which confirms a conjecture of Z.-W. Sun. We also establish the following new congruences:
\begin{align*}\sum_{k=0}^{(p-1)/2}\frac{\binom{2k}k\binom{3k}k}{27^k}\equiv&\l(\frac p3\r)\frac{2^p+1}3\pmod{p^2},
\\\sum_{k=0}^{(p-1)/2}\frac{\binom{6k}{3k}\binom{3k}k}{(2k+1)432^k}\equiv&\l(\frac p3\r)\frac{3^p+1}4\pmod{p^2},
\\\sum_{k=0}^{(p-1)/2}\frac{\binom{4k}{2k}\binom{2k}k}{(2k+1)64^k}\equiv&\l(\frac{-1}p\r)2^{p-1}\pmod{p^2}.
\end{align*}
Note that in 2003 Rodriguez-Villeguez posed conjectures on
$$\sum_{k=0}^{p-1}\frac{\binom{2k}k^2}{16^k},\ \sum_{k=0}^{p-1}\frac{\binom{2k}k\binom{3k}k}{27^k},\ \sum_{k=1}^{p-1}\frac{\binom{4k}{2k}\binom{2k}k}{64^k},\ \sum_{k=1}^{p-1}\frac{\binom{6k}{3k}\binom{3k}k}{432^k}$$
modulo $p^2$ which were later proved.
\end{abstract}
\maketitle

\section{Introduction}
\setcounter{lemma}{0}
\setcounter{theorem}{0}
\setcounter{corollary}{0}
\setcounter{remark}{0}
\setcounter{equation}{0}

Let $p>3$ be a prime. In 2003, via his analysis
of the p-adic analogues of Gaussian hypergeometric series and the Calabi-
Yau manifolds, Rodriguez-Villegas \cite{RV} conjectured the following congruences:
\begin{gather*}\sum_{k=0}^{p-1}\f{\bi{2k}k^2}{16^k}\eq\l(\f{-1}p\r)\pmod{p^2},\ \sum_{k=0}^{p-1}\frac{\binom{2k}k\binom{3k}k}{27^k}\equiv\l(\frac p3\r)\pmod{p^2},
\\\sum_{k=0}^{p-1}\frac{\binom{4k}{2k}\binom{2k}k}{64^k}\equiv\l(\frac{-2}p\r)\pmod{p^2},
\ \sum_{k=0}^{p-1}\frac{\binom{6k}{3k}\binom{3k}k}{432^k}\equiv\l(\frac {-1}p\r)\pmod{p^2},
\end{gather*}
where $(\f{\cdot}p)$ denotes the Jacobi symbol.
They were soon proved by E. Mortenson \cite{M1, M2} via the Gross-Koblitz formula and the $p$-adic $\Gamma$-function. Note that
\begin{gather*}\bi{-1/2}k^2=\f{\bi{2k}k^2}{16^k},\ \ \bi{-1/3}k\bi{-2/3}k=\f{\bi{2k}k\bi{3k}k}{27^k},
\\\bi{-1/4}k\bi{-3/4}k=\f{\bi{4k}{2k}\bi{2k}k}{64^k},\ \ \bi{-1/6}k\bi{-5/6}k=\f{\bi{6k}{3k}\bi{3k}k}{432^k}
\end{gather*}
for all $k\in\N=\{0,1,2,\ldots\}$.
In 2011 Z. W. Sun [Su11] showed further that
$$\sum_{k=1}^{p-1}\f{\bi{2k}k^2}{16^k}\eq\l(\f{-1}p\r)-p^2E_{p-3}\pmod{p^3}$$
and
\begin{equation}\label{1.1}\sum_{k=1}^{(p-1)/2}\f{\bi{2k}k^2}{16^k}\eq\l(\f{-1}p\r)+p^2E_{p-3}\pmod{p^3},
\end{equation}
where $E_0,E_1,E_2,\ldots$ are the Euler numbers given by 
$$E_0=1,\ \t{and}\ E_n=-\sum_{k=1}^{\lfloor n/2\rfloor}\binom{n}{2k}E_{n-2k}\ (n=1,2,3,\ldots).$$
He also conjectured that
\begin{align*}\sum_{k=0}^{p-1}\f{\bi{2k}k\bi{3k}k}{(2k+1)27^k}\eq&\l(\f p3\r)\pmod{p^2},
\\\sum_{k=0}^{p-1}\f{\bi{4k}{2k}\bi{2k}k}{(2k+1)64^k}\eq&\l(\f{-1}p\r)-3p^2E_{p-3}\pmod{p^3},
\\\sum_{k=0}^{p-1}\f{\bi{6k}{3k}\bi{3k}k}{(2k+1)432^k}\eq&\l(\f p3\r)\pmod{p^2},
\end{align*}
which were confirmed by Z.-H. Sun \cite{S16}. Note that Z.-W. Sun \cite{Su14} determined
$$\sum_{k=0}^{(p-3)/2}\f{\bi{2k}k^2}{(2k+1)16^k}\ \ \t{and}\ \  \sum_{k=(p+1)/2}^{p-1}\f{\bi{2k}k^2}{(2k+1)16^k}$$
modulo $p^3$.

In this paper we first establish the following result.

\begin{theorem}\label{Th1.1} Let $p$ be any odd prime.

{\rm (i)} We have
\begin{equation}\label{1.2}\sum_{k=0}^{\lfloor 3p/4\rfloor}\frac{\binom{2k}k^2}{16^k}\equiv\begin{cases}1\pmod{p^3}&\t{if}\ p\eq1\pmod 4,
\\-1+p^2/(2\bi{(p-3)/2}{(p-3)/4}^2)\pmod{p^3}&\t{if}\ p\eq3\pmod 4.\end{cases}
\end{equation}

{\rm (ii)} For each $a=2,3,4,\ldots$, we have
\begin{equation}\label{1.3}\sum_{k=0}^{\lfloor\frac34p^a\rfloor}\frac{\binom{2k}k^2}{16^k}\equiv\l(\frac{-1}{p^a}\r)\pmod{p^3}.
\end{equation}
\end{theorem}
\begin{remark}\label{Rem1.1}. Part (i) in the case $p\eq1\pmod 4$ and part (ii) were conjectured by Sun \cite{Su11}.
\end{remark}

Our second theorem is as follows.
\begin{theorem}\label{Th1.2} Let $p>3$ be a prime. Then we have
\begin{equation}\label{1.4}
\sum_{k=0}^{(p-1)/2}\frac{\binom{2k}k\binom{3k}k}{27^k}\equiv\l(\frac p3\r)\frac{2^p+1}3\pmod{p^2},
\end{equation}
\begin{equation}\label{1.5}
\sum_{k=0}^{(p-1)/2}\frac{\binom{6k}{3k}\binom{3k}k}{(2k+1)432^k}\equiv\l(\frac p3\r)\frac{3^p+1}4\pmod{p^2},
\end{equation}
\begin{equation}\label{1.6}
\sum_{k=0}^{(p-1)/2}\frac{\binom{4k}{2k}\binom{2k}k}{(2k+1)64^k}\equiv\l(\frac{-1}p\r)2^{p-1}\pmod{p^2}.
\end{equation}
\end{theorem}
\begin{remark}\label{Rem1.2} We are also able to show the congruence
$$\sum_{k=0}^{(p-1)/2}\frac{\binom{2k}k\binom{3k}k}{(2k+1)27^k}\equiv\l(\frac p3\r)(3^p+2-2^{p+1})\pmod {p^2}$$
for any prime $p>3$.
\end{remark}
\section{Proof of Theorem 1.1}
\setcounter{lemma}{0}
\setcounter{theorem}{0}
\setcounter{corollary}{0}
\setcounter{remark}{0}
\setcounter{equation}{0}

\begin{lemma}\label{Lem2.1} {\rm (Sun \cite[(1.4)]{Su11})} For any prime $p>3$ we have
\begin{equation}\label{2.1}\sum_{k=1}^{(p-1)/2}\f{4^k}{k^2\bi{2k}k}\eq (-1)^{(p-1)/2}4E_{p-3}\pmod p.
\end{equation}
\end{lemma}

\medskip
\noindent{\it Proof of Theorem} 1.1(i). In view of (\ref{1.1}), (\ref{1.2})
has the following equation form:
\begin{equation}\label{2.2}\sum_{k=(p+1)/2}^{\lfloor 3p/4\rfloor}\f{\bi{2k}k^2}{16^k}\eq-p^2E_{p-3}+\f{1-(-1)^{(p-1)/2}}2\cdot\f{p^2}{2\bi{(p-3)/2}{\lfloor p/4\rfloor}^2}\pmod{p^3}.
\end{equation}
By \cite[Lemma 2.1]{Su11},
$$k\binom{2k}k\binom{2(p-k)}{p-k}\equiv(-1)^{\lfloor2k/p\rfloor-1}2p\pmod{p^2}\ \t{for all}\ k=1,\ldots,p-1.$$
Thus
\begin{align*}\sum_{k=(p+1)/2}^{\lfloor 3p/4\rfloor}\frac{\binom{2k}k^2}{16^k}
\equiv&\sum_{k=(p+1)/2}^{\lfloor 3p/4\rfloor}\frac{4p^2}{k^2\binom{2(p-k)}{p-k}^216^k}
=\sum_{j=\lfloor p/4\rfloor+1}^{(p-1)/2}\frac{4p^2}{(p-j)^2\binom{2j}{j}^216^{p-j}}
\\\eq&\frac{p^2}4\sum_{j=\lfloor p/4\rfloor+1}^{(p-1)/2}\frac{16^j}{j^2\binom{2j}j^2}\pmod{p^3}
\end{align*}
\and hence we have reduced (\ref{2.2}) to the following simpler form
\begin{equation}\label{2.3} \sum_{k=\lfloor n/2\rfloor+1}^n\f{16^k}{k^2\bi{2k}k^2}\eq-4E_{p-3}+\f{1-(-1)^{n}}{\bi{n-1}{\lfloor n/2\rfloor}^2}\pmod p,\end{equation}
where $n=(p-1)/2$.

For each $k=0,\ldots,n$, clearly
$$\bi nk\eq\bi{-1/2}k=\f{\bi{2k}k}{(-4)^k}\pmod{p}.$$
Thus
$$\sum_{k=\lfloor n/2\rfloor+1}^{n}\frac{16^k}{k^2\binom{2k}k^2}\equiv\sum_{k=\lfloor n/2\rfloor+1}^{n}\frac1{k^2\binom nk^2}\equiv4\sum_{k=\lfloor n/2\rfloor+1}^{n}\frac1{\binom{n-1}{k-1}^2}\pmod p.$$
Note that $$\sum_{k=\lfloor n/2\rfloor+1}^{n}\frac1{\binom{n-1}{k-1}^2}=\sum_{k=\lfloor n/2\rfloor}^{n-1}\frac1{\binom{n-1}{k}^2}=\frac12\sum_{k=0}^{n-1}\frac1{\binom{n-1}{k}^2}+\f{1-(-1)^n}{4\bi{n-1}{\lfloor n/2\rfloor}^2}$$
and
\begin{equation}\label{2.4}
\sum_{k=0}^{n-1}\f1{\bi{n-1}k^2}=\f{2n^2}{n+1}\sum_{k=1}^{n}\f1{k\bi{2n+1-k}{n-k}}
\end{equation}
(cf. \cite{SWZ}). So we have
\begin{align*}&\sum_{k=\lfloor n/2\rfloor+1}^{n}\frac{16^k}{k^2\binom{2k}k^2}-\f{1-(-1)^n}{\bi{n-1}{\lfloor n/2\rfloor}^2}
\\\equiv&\frac{4n^2}{n+1}\sum_{k=1}^{n}\frac1{k\binom{2n+1-k}{n-k}}\equiv2\sum_{k=1}^{n}\frac{1}{k\binom{-k}{n-k}}\pmod p.
\end{align*}
Observe that
\begin{align*}
\sum_{k=1}^n\f1{k\bi{-k}{n-k}}=&\sum_{k=1}^n\f{(-1)^{n-k}}{k\bi{n-1}{k-1}}=n\sum_{k=1}^n\f{(-1)^{n-k}}{k^2\bi{n}k}
\\\eq&\f{(-1)^{n-1}}2\sum_{k=1}^n\f{4^k}{k^2\bi{2k}k}\pmod p.
\end{align*}
Therefore, with the help of Lemma \ref{Lem2.1}, we finally obtain
$$ \sum_{k=\lfloor n/2\rfloor+1}^n\f{16^k}{k^2\bi{2k}k^2}-\f{1-(-1)^n}{\bi{n-1}{\lfloor n/2\rfloor}^2}
\eq(-1)^{n-1}\sum_{k=1}^n\f{4^k}{k^2\bi{2k}k}\eq-4E_{p-3}\pmod p.$$
This proves (\ref{2.3}) and hence (\ref{1.2}) follows. \qed

Now we give a lemma which is a natural extension of (\ref{1.1}).

\begin{lemma}\label{Lem2.3} Let $p>3$ be a prime and let $a$ be any positive integer. Then
\begin{equation}\label{2.5}\sum_{k=0}^{(p^a-1)/2}\frac{\binom{2k}k^2}{16^k}\equiv\l(\frac{-1}{p^a}\r)+\l(\f{-1}{p^{a-1}}\r)p^2E_{p-3}\pmod{p^3}.
\end{equation}
\end{lemma}
\Proof. Theorem 1.2 of Sun \cite{Su13} states that for any $d=0,\ldots,(p-1)/2$ we have
$$\sum_{k=0}^{(p-1)/2}\f{\bi{2k}k\bi{2k}{k+d}}{16^k}\eq\l(\f{-1}p\r)+\f{(-1)^d}4p^2E_{p-3}\l(d+\f12\r)\pmod{p^3},$$
where $E_{p-3}(x)$ denotes the Euler polynomial of the degree $p-3$.

In the case $d=0$ this yields (1.1). Modifying this proof of (\ref{1.1}) slightly we immediately get (\ref{2.5}). \qed
\smallskip

In 1852, Kummer proved that for any $m,n\in\N$ the $p$-adic valuation of the binomial coefficient $\binom{m+n}m$ is equal to the number of {\it carry-overs} when performing the addition of $m$ and $n$ written in base $p$.

\begin{lemma}\label{Lem2.4} Let $p$ be an odd prime and let $a\in\mathbb{Z}^{+}$. For any $k=1,2,\ldots,(p^a-1)/2$, we have
$$\ord_p\binom{p^a-k}{\frac{p^a-1}2-k}\leq a-1.$$
\end{lemma}
{\it Proof}. It is well known that
$$\ord_p(n!)=\sum_{j=1}^\infty\l\lfloor\f{n}{p^j}\r\rfloor.$$
Thus
$$\ord_p\binom{p^a-k}{\frac{p^a-1}2-k}=\sum_{j=1}^{a-1}\l(\l\lfloor\f{p^a-k}{p^j}\r\rfloor-\l\lfloor\f{(p^a+1)/2}{p^j}\r\rfloor-\l\lfloor\f{(p^a-1)/2-k}{p^j}\r\rfloor\r)$$
does not exceed $a-1$ as each term in the sum is at most one.
This concludes the proof. \qed

\medskip
\noindent{\it Proof of Theorem} 1.1(ii). In view of  Lemma \ref{Lem2.3}, we just need to verify that
\begin{equation}\label{2.6}
\sum_{k=(p^a+1)/2}^{\lfloor3p^a/4\rfloor}\frac{\binom{2k}k^2}{16^k}\equiv\l(\f{-1}{p^{a-1}}\r)p^2E_{p-3}\pmod{p^3}.
\end{equation}
Let $k$ and $l$ be positive integers with $k+l=p^a$ and $0<l<p^a/2$. Then $$\frac{\binom{2k}k^2}{\binom{2p^a-2}{p^a-1}^2}=\frac{(2p^a-2l)!^2}{(2p^a-2)!^2}\l(\frac{(p^a-1)!}{(p^a-l)!}\r)^4=\frac{\prod_{0<i<l}(p^a-i)^4}{\prod_{1<j<2l}(2p^a-j)^2}$$ and hence$$\frac{\binom{2k}k^2}{\binom{2p^a-2}{p^a-1}^2}\cdot\frac{(2l-1)!^2}{(l-1)!^4}=\frac{\prod_{0<i<l}(1-p^a/i)^4}{\prod_{1<j<2l}(1-2p^a/j)^2}\equiv1\pmod{p}.$$
Note that $$\binom{2p^a-2}{p^a-1}^2=p^{2a}\prod_{j=2}^{p^a-1}\l(\frac{2p^a-j}j\r)^2\equiv p^{2a}\pmod{p^{2a+1}}$$
and $$\binom{2k}k^2=\binom{p^a+(2k-p^a)}{0p^a+k}^2\equiv\binom{2k-p^a}k^2=0\pmod{p^2}$$
 by Lucas' theorem. So we have $$\frac{l^2}4\binom{2l}l^2=\frac{(2l-1)!^2}{(l-1)!^4}\not\equiv0\pmod{p^{2a}}$$ and $$\binom{2k}k^2\equiv p^{2a}\frac{(l-1)!^4}{(2l-1)!^2}=\frac{4p^{2a}}{l^2\binom{2l}l^2}\pmod{p^3}.$$
Therefore
\begin{align*}
\sum_{k=(p^a+1)/2}^{\lfloor3p^a/4\rfloor}\frac{\binom{2k}k^2}{16^k}\equiv&\sum_{k=(p^a+1)/2}^{\lfloor3p^a/4\rfloor}\frac{4p^{2a}}{16^k(p^a-k)^2\binom{2(p^a-k)}{p^a-k}^2}
\\\eq&\frac{p^{2a}}4\sum_{l=\lfloor p^a/4\rfloor+1}^{(p^a-1)/2}\frac{16^l}{l^2\binom{2l}l^2}\pmod{p^3}.
\end{align*}
For $k=1,\ldots,(p^a-1)/2$, clearly
\begin{align*}
\frac{\binom{(p^a-1)/2}{k}}{\binom{2k}k/{(-4)^k}}=&\frac{\binom{(p^a-1)/2}{k}}{\binom{-1/2}k}=\prod_{j=0}^{k-1}\frac{(p^a-1)/2-j}{-1/2-j}
\\=&\prod_{j=0}^{k-1}\l(1-\frac{p^a}{2j+1}\r)\equiv1\pmod p.
\end{align*}
Thus
\begin{align*}
\sum_{k=(p^a+1)/2}^{\lfloor3p^a/4\rfloor}\frac{\binom{2k}k^2}{16^k}\equiv&\frac{p^{2a}}4\sum_{k=\lfloor p^a/4\rfloor+1}^{(p^a-1)/2}\frac{1}{k^2\binom{(p^a-1)/2}k^2}
\\\equiv& p^{2a}\sum_{k=\lfloor p^a/4\rfloor+1}^{(p^a-1)/2}\frac{1}{\binom{(p^a-3)/2}{k-1}^2}\pmod{p^3}.
\end{align*}
So (\ref{2.6}) is reduced to
\begin{equation}\label{2.7}
p^{2a-2}\sum_{k=\lfloor p^a/4\rfloor}^{(p^a-3)/2}\frac1{\binom{(p^a-3)/2}{k}^2}\equiv-\l(\f{-1}{p^{a-1}}\r)E_{p-3}\pmod p.
\end{equation}
If $p^a\equiv1\pmod 4$, then $(p^a-3)/2$ is odd and hence $$\sum_{k=\lfloor p^a/4\rfloor}^{(p^a-3)/2}\frac1{\binom{(p^a-3)/2}{k}^2}=\frac12\sum_{k=0}^{(p^a-3)/2}\frac1{\binom{(p^a-3)/2}{k}^2}.$$
If $p^a\equiv3\pmod 4$, then $a\in\{3,5,\ldots\}$ and $$\sum_{k=\lfloor p^a/4\rfloor}^{(p^a-3)/2}\frac1{\binom{(p^a-3)/2}{k}^2}=\frac12\sum_{k=0}^{(p^a-3)/2}\frac1{\binom{(p^a-3)/2}{k}^2}+\frac12\cdot\frac1{\binom{(p^a-3)/2}{(p^a-3)/4}^2}.$$
In the case $p^a\eq3\pmod4$, as the fractional parts of $(p^a-3)/(2p)$ and $(p^a-3)/(4p)$ are $(p-3)/(2p)$ and $(p-3)/(4p)$ respectively, we have $$\l\lfloor\frac{(p^a-3)/2}p\r\rfloor=2\l\lfloor\frac{(p^a-3)/4}p\r\rfloor$$ and hence $$\ord_p\binom{(p^a-3)/2}{(p^a-3)/4}^2=2\sum_{j=1}^{a-1}\l(\l\lfloor\frac{(p^a-3)/2}{p^j}\r\rfloor-2\l\lfloor\frac{(p^a-3)/4}{p^j}\r\rfloor\r)<2a-2.$$
No matter $p^a\eq1\pmod 4$ or not, we always have
$$p^{2a-2}\sum_{k=\lfloor p^a/4\rfloor}^{(p^a-3)/2}\frac1{\binom{(p^a-3)/2}{k}^2}\equiv\frac{p^{2a-2}}2\sum_{k=0}^{(p^a-3)/2}\frac1{\binom{(p^a-3)/2}{k}^2}\pmod{p}.$$
So (\ref{2.7}) has the following equivalent form:
\begin{equation}\label{2.8}p^{2a-2}\sum_{k=0}^{(p^a-3)/2}\frac1{\binom{(p^a-3)/2}{k}^2}\equiv-2\l(\f{-1}{p^{a-1}}\r)E_{p-3}\pmod p.\end{equation}

The identity (\ref{2.4}) with $n=(p^a-1)/2$ yields that $$\sum_{k=0}^{(p^a-3)/2}\frac1{\binom{(p^a-3)/2}{k}^2}=\frac{2((p^a-1)/2)^2}{(p^a+1)/2}\sum_{k=1}^{(p^a-1)/2}\frac1{k\binom{p^a-k}{(p^a-1)/2-k}}.$$
So (\ref{2.8}) is reduced to
\begin{equation}\label{2.9}
p^{2a-2}\sum_{k=1}^{(p^a-1)/2}\frac1{k\binom{p^a-k}{(p^a+1)/2}}\equiv-2\l(\f{-1}{p^{a-1}}\r)E_{p-3}\pmod p.
\end{equation}
In view of Lemma \ref{Lem2.4}, if $1\ls k\ls (p^a-1)/2$ and $p^{a-1}\nmid k$, then $$\frac{p^{2a-2}}{k\binom{p^a-k}{(p^a+1)/2}}\equiv0\pmod p.$$
Thus
\begin{align*}
p^{2a-2}\sum_{k=1}^{(p^a-1)/2}\frac1{k\binom{p^a-k}{(p^a+1)/2}}\equiv& p^{2a-2}\sum_{j=1}^{(p-1)/2}\frac1{p^{a-1}j\binom{p^a-p^{a-1}j}{(p^a+1)/2}}
\\=&\frac{p^a+1}2\sum_{j=1}^{(p-1)/2}\frac1{j(p-j)\binom{p^a-p^{a-1}j-1}{(p^a-1)/2}}
\\\equiv&-\frac12\sum_{j=1}^{(p-1)/2}\frac1{j^2\binom{p^a-p^{a-1}j-1}{(p^a-1)/2}}\pmod{p}.
\end{align*}
For each $j=1,\ldots,(p-1)/2$, by Lucas' theorem we have
\begin{align*}
\binom{p^{a-1}(p-j)-1}{(p^a-1)/2}= &\bi{p^{a-1}(p-1-j)+p^{a-1}-1}{p^{a-1}(p-1)/2+(p^{a-1}-1)/2}
\\\eq&\bi{p-1-j}{(p-1)/2}\bi{p^{a-1}-1}{(p^a-1)/2}
\\\eq&(-1)^{(p^{a-1}-1)/2}\binom{p-j-1}{(p-1)/2}\pmod p,
\end{align*}
also
\begin{align*}\bi{p-j-1}{(p-1)/2}=&\bi{p-1-j}{(p-1)/2-j}=(-1)^{(p-1)/2-j}\bi{-p+(p-1)/2}{(p-1)/2-j}
\\\eq&(-1)^{(p-1)/2-j}\bi{(p-1)/2}j\eq(-1)^{(p-1)/2-j}\bi{-1/2}j
\\=&(-1)^{(p-1)/2}\f{\bi{2j}j}{4^j}\pmod p.\end{align*}
Therefore
\begin{align*}
&p^{2a-2}\sum_{k=1}^{(p^a-1)/2}\frac1{k\binom{p^a-k}{(p^a+1)/2}}\equiv\frac{(-1)^{(p^{a-1}+1)/2}}2\sum_{j=1}^{(p-1)/2}\frac1{j^2\binom{p-j-1}{(p-1)/2}}
\\\equiv&\frac{(-1)^{(p^{a-1}+1)/2}}2(-1)^{(p-1)/2}\sum_{j=1}^{(p-1)/2}\frac{4^j}{j^2\binom{2j}j}\pmod p.
\end{align*}
This, together with (\ref{2.1}), yields the desired (\ref{2.9}).

The proof of Theorem 1.1(ii) is now complete. \qed

\section{Proof of Theorem 1.2}
\setcounter{lemma}{0}
\setcounter{theorem}{0}
\setcounter{corollary}{0}
\setcounter{remark}{0}
\setcounter{equation}{0}
\setcounter{conjecture}{0}
\begin{lemma}\label{Lem3.1} Let $p>3$ be a prime, and $m\in\{1,2,\ldots,(p-1)/2\}$. For any $p$-adic integer $t$, we have
\begin{align}\label{3.1}
 \binom{m+pt-1}{(p-1)/2}\binom{-1-pt-m}{(p-1)/2}\equiv\frac{pt}m\pmod{p^2}.
 \end{align}
\end{lemma}
\Proof. Since
\begin{align*}
\binom{m+pt-1}{(p-1)/2}
=&\frac{\prod_{r=0}^{m-1}(pt+r)\times\prod_{s=1}^{(p-1)/2-m}(pt-s)}{((p-1)/2)!}
\\\eq&\frac{(m-1)!pt(-1)^{(p-1)/2-m}((p-1)/2-m)!}{((p-1)/2)!}\pmod{p^2},
\end{align*}
and
\begin{align*}
&\binom{-m-pt-1}{(p-1)/2}=\frac{\prod_{j=1}^{(p-1)/2}(-m-pt-j)}{((p-1)/2)!}
\\\eq&\frac{(-1)^{(p-1)/2}(m+1)(m+2)\cdots(m+(p-1)/2)}{((p-1)/2)!}\pmod{p},
\end{align*}
we have
\begin{align*}
&\binom{m+pt-1}{(p-1)/2}\binom{-m-pt-1}{(p-1)/2}
\\\eq&\frac{pt(m-1)!(-1)^m((p-1)/2-m)!(m+1)(m+2)\cdots(m+(p-1)/2)}{((p-1)/2)!^2}
\\=&\frac{pt}m\frac{(-1)^m((p-1)/2-m)!(m+(p-1)/2)!}{((p-1)/2)!^2}=\frac{pt}m(-1)^m\frac{\binom{p-1}{(p-1)/2}}{\binom{p-1}{(p-1)/2+m}}
\\\eq&\frac{pt}m(-1)^m(-1)^{(p-1)/2}(-1)^{(p-1)/2+m}=\frac{pt}m\pmod{p^2}.
\end{align*}
This concludes the proof. \qed

\begin{remark}\label{Rem3.1} Let $p>3$ be a prime and $m\in\{(p+1)/2,\ldots,p-1\}$. For any $p$-adic integer $t$, by Lemma \ref{Lem3.1} we have
 \begin{align*}&\binom{m+pt-1}{(p-1)/2}\binom{-1-pt-m}{(p-1)/2}
 \\=&\binom{(m-p)+p(t+1)-1}{(p-1)/2}\binom{-1-p(t+1)-(m-p)}{(p-1)/2}
 \\\equiv&\frac{p(t+1)}{m-p}\eq\f{p(t+1)}m\pmod{p^2}.\end{align*}
\end{remark}

\begin{lemma}\label{Lem3.2} Let $p>3$ be a prime. For $k\in\{1,2,\ldots,p-1\}$ and $p$-adic integer $t$, we have
\begin{equation}\label{tk}
 \binom{pt}{k}\binom{-1-pt}{k}\equiv-\frac{p^2t^2}{k^2}-\frac{pt}k\pmod{p^3}.
 \end{equation}
\end{lemma}
\Proof. This is almost trivial. In fact,
\begin{align*}
\binom{pt}k\binom{-1-pt}k=&\frac{pt}{pt-k}\bi{-1+pt}k\bi{-1-pt}k
\\\eq& \f{pt}{pt-k}\bi{-1}k^2=\f{pt(p^2t^2+ptk+k^2)}{(pt)^3-k^3}
\\\eq&-\frac{p^2t^2}{k^2}-\frac{pt}k\pmod{p^3}.
\end{align*}
This proves (\ref{tk}).\qed
\medskip

Recall that those $H_n=\sum_{0<k\ls n}1/k$ with $n\in\N$ are called harmonic numbers. If a prime $p$ does not divide an integer $a$, then we let $q_p(a)$ denote the
Fermat quotient $(a^{p-1}-1)/p$.

\begin{lemma}\label{Lehmer} {\rm (Lemma \cite{L})}. For any prime $p>3$, we have
\begin{gather*}H_{\lfloor p/2\rfloor}\eq-2q_p(2)\pmod p,\  H_{\lfloor p/4\rfloor}\eq-3q_p(2)\pmod p,
\\H_{\lfloor p/3\rfloor}\eq-\f32q_p(3)\pmod p \ \t{and}\ H_{\lfloor p/6\rfloor}\eq-2q_p(2)-\f 32q_p(3)\pmod p,
\end{gather*}
where $q_p(2)=(2^{p-1}-1)/p$ and $q_p(3)=(3^{p-1}-1)/p$ stand for the Fermat quotients.
\end{lemma}

For $n\in\N$, define 
$$S_n(x)=\sum_{k=0}^n\binom{x}k\binom{-1-x}k\ \ \mbox{and}\ \ T_n(x)=\sum_{k=0}^n\binom{x}k\binom{-1-x}k\frac{1+2x}{1+2k}.$$
By \cite[(2.2)]{S16} with $a=x+1$ and $b=0$, we have
\begin{align}\label{3.3}
S_n(x)+S_n(x+1)=2\binom{x}n\binom{-2-x}{n}.
\end{align}
By \cite[(2.2)]{S16} with $b=2$, we get
\begin{align}\label{3.5}
T_n(x)-T_n(x-1)=2\binom{x-1}n\binom{-x-1}n.
\end{align}

\medskip
\noindent{\it Proof of Theorem 1.2}. For any $p$-adic integer $a$, we let $\langle a\rangle_p$
denote the least nonnegative integer $r$ with $a\eq r\pmod p$. For convenience, we also set $n=(p-1)/2$. 

(i)  For any $p$-adic integer $a\not\equiv0\pmod p$,  by using (\ref{3.3}) we get
\begin{align*}&S_n(a)-(-1)^{\langle a\rangle_p}S_n(a-\langle a\rangle_p)
\\=&\sum_{k=0}^{\langle a\rangle_p-1}(-1)^k(S_n(a-k)+S_n(a-k-1))
\\=&\sum_{k=0}^{\langle a\rangle_p-1}(-1)^k 2\binom{a-k-1}n\binom{k-a-1}{n}
\end{align*}
and hence
\begin{align*}&S_n(a)-(-1)^{\langle a\rangle_p}S_n(pt)
\\=&2\sum_{k=0}^{\langle a\rangle_p-1}(-1)^k\bi{\langle a\rangle_p+pt-k-1}{n}\bi{-1-pt-(\langle a\rangle_p-k)}{n},
\end{align*}
where $t:=(a-\langle a\rangle_p)/p$. By Lemma \ref{Lem3.2},
\begin{align*}\label{3.4}
S_n(pt)=\sum_{k=0}^n\binom{pt}k\binom{-1-pt}k\equiv1-\sum_{k=1}^n\frac{pt}k=1-ptH_n\pmod{p^2}.
\end{align*}
So, with helps of Lemma 3.1 and Remark 3.1,  we have
\begin{equation}\label{delta}S_n(a)-(-1)^{\langle a\rangle_p}(1-ptH_n)
\eq2\sum_{k=0}^{\langle a\rangle_p-1}(-1)^k\f{p(t+\da_k)}{\langle a\rangle_p-k}\pmod{p^2},
\end{equation}
where $\da_k$ takes $1$ or $0$ according as $\langle a\rangle_p-k>p/2$ or not.

Observe that
$$\sum_{k=0}^{(p-1)/2}\f{\bi{2k}k\bi{3k}k}{27^k}=\sum_{k=0}^n\bi{-1/3}k\bi{-2/3}k=S_n(a)$$
with $a=-1/3$. Note that
$$\langle a\rangle_p=\begin{cases}(p-1)/3&\t{if}\ p\eq1\pmod 3,
\\(2p-1)/3&\t{if}\ p\eq2\pmod 3.\end{cases}$$
Hence
$$t:=\f{a-\langle a\rangle_p}p=\begin{cases}-1/3&\t{if}\ p\eq 1\pmod 3,\\-2/3&\t{if}\ p\eq2\pmod 3.\end{cases}$$

{\it Case} 1. $p\equiv1\pmod3$.

In this case, $\langle a\rangle_p=(p-1)/3$,  $t=-1/3$, and $\da_k=0$ for all $k=0,\ldots,\langle a\rangle_p-1$. So
we have
\begin{align*}&S_n\l(-\f13\r)-(-1)^{(p-1)/3}(1-ptH_n)
\\\eq& 2pt(-1)^{(p-1)/3}\sum_{j=1}^{(p-1)/3}\f{(-1)^j}j=2pt\l(H_{(p-1)/6}-H_{(p-1)/3}\r)\pmod{p^2}.
\end{align*}
Combining this with Lemma \ref{Lehmer} and recalling that $t=-1/3$, we immediately obtain the desired congruence
$$S_n\l(-\f13\r)\eq1+\f23 p\,q_p(2)\pmod{p^2}.$$

{\it Case} 2. $p\equiv2\pmod 3$.

In this case, we have $\langle a\rangle_p=(2p-1)/3$, $t=-2/3$ and $$\da_k=\begin{cases} 1&\t{if}\  0\leq k<(p+1)/6,\\ 0&\t{if}\ (p+1)/6\leq k\leq\langle a\rangle_p-1.\end{cases}$$ So we have
\begin{align*}
&S_n\l(-\f13\r)-(-1)^{(2p-1)/3}(1-ptH_n)
\\\eq&2p(t+1)\sum_{k=0}^{(p-5)/6}\f{(-1)^k}{\langle a\rangle_p-k}+2pt\sum_{k=(p+1)/6}^{(2p-4)/3}\frac{(-1)^k}{\langle a\rangle_p-k}
\\=& 2p(t+1)(-1)^{(2p-1)/3}\sum_{j=(p+1)/2}^{(2p-1)/3}\f{(-1)^j}j+2pt(-1)^{(2p-1)/3}\sum_{j=1}^{(p-1)/2}\frac{(-1)^j}j
\\=&-2p(t+1)\sum_{j=1}^{(2p-1)/3}\f{(-1)^j}j+2p\sum_{j=1}^{(p-1)/2}\frac{(-1)^j}j
\\=&-2p(t+1)\l(H_{\lfloor p/3\rfloor}-H_{\lfloor {2p}/3\rfloor}\r)+2p\l(H_{\lfloor p/4\rfloor}-H_{\lfloor p/2\rfloor}\r)\pmod{p^2}.
\end{align*}
Note that
$$H_{\lfloor {2p}/3\rfloor}=\sum_{k=1}^{(p-1)/2}\l(\f1k+\f1{p-k}\r)-\sum_{j=1}^{(p-1)/3}\f1{p-j}\eq H_{\lfloor p/3\rfloor}\pmod p.$$
Therefore,
$$S_n\l(-\f13\r)+1-ptH_{\lfloor p/2\rfloor}\eq 2p\l(H_{\lfloor p/4\rfloor}-H_{\lfloor p/2\rfloor}\r)\pmod{p^2}.$$
 This, together with Lemma \ref{Lehmer} and the fact that $t=-2/3$, yields the desired congruence $$S_n\l(-\f13\r)\eq-1-\f23 p\,q_p(2)\pmod{p^2}.$$

In view of the above, we have completed the proof of (\ref{1.4}).

(ii) For any $p$-adic integer $a$ with $a(2a+1)\not\eq0\pmod p$, if we set $t=(a-\langle a\rangle_p)/p$ then
by (\ref{3.5}) we have
\begin{align*}T_n(a)-T_n(pt)=&\sum_{k=1}^{\langle a\rangle_p}(T_n(a-k+1)-T_n(a-k))
\\=&\sum_{k=1}^{\langle a\rangle_p}2\bi{a-k}n\bi{k-a-2}n
\\=&2\sum_{k=1}^{\langle a\rangle_p}\bi{m_k+pt-1}n\bi{-1-pt-m_k}n,
\end{align*}
where $m_k=\langle a\rangle_p-k+1$.
In view of Lemmas \ref{Lem3.2} and \ref{Lehmer},
\begin{align*}
T_n(pt)-(1+2pt)=&\sum_{k=0}^n\binom{pt}k\binom{-1-pt}k\frac{1+2pt}{1+2k}-(1+2pt)
\\\eq&\binom{pt}n\binom{-1-pt}n\frac{1+2pt}p-\sum_{k=1}^{n-1}\frac{pt}{k(1+2k)}
\\\eq&\l(-\frac{p^2t^2}{n^2}-\frac{pt}n\r)\frac{1+2pt}p-\sum_{k=1}^{n-1}\frac{pt}{k(1+2k)}
\\\eq&2t+2pt-pt\sum_{k=1}^{n-1}\frac1k+2pt\sum_{k=1}^{n-1}\frac1{2k+1}
\\\eq&2t-2pt-ptH_n+2pt\l(H_{p-1}-\f{H_n}2\r)
\\\eq&2t-2pt+4ptq_p(2)\pmod{p^2}
\end{align*}
and hence
$$T_n(pt)\equiv1+2t+4ptq_p(2)\pmod{p^2}.$$
Therefore, with the helps of Lemma \ref{Lem3.1} and Remark 3.1, we have
\begin{align*}&T_n(a)-(1+2t+4ptq_p(2))
\\\eq&2\sum_{k=1}^{\langle a\rangle_p}\bi{m_k+pt-1}n\bi{-1-pt-m_k}n
\\\eq&2\sum_{k=1}^{\langle a\rangle_p}\f{p(t+\da_k)}{m_k}
=2\sum_{j=1}^{\langle a\rangle_p}\f{pt}j+2\sum^{\langle a\rangle_p}_{j=1\atop j>p/2}\f1j\pmod {p^2},
\end{align*}
where $\da_k$ takes $1$ or $0$ according as $m_k>p/2$ or not.
Below we deal with $a=-1/6,-1/4$.

Clearly,
$$H_{p-k}=H_{p-1}-\sum_{0<j<k}\f1{p-j}\eq H_{k-1}\pmod p$$
for all $k=1,\ldots,p-1$. Thus, with the help of Lemma \ref{Lehmer} we have
$$H_{\lfloor{3p}/4\rfloor}\equiv H_{p-1-\lfloor 3p/4\rfloor}=H_{\lfloor p/4\rfloor}\eq-3q_p(2)\pmod p$$
and
$$H_{\lfloor {5p}/6\rfloor}\equiv H_{p-1-\lfloor 5p/6\rfloor}=H_{\lfloor p/6\rfloor}\eq-2q_p(2)-\frac32q_p(3)\pmod p.$$

{\it Case} I. $\langle a\rangle_p<n$.

If $a=-1/6$, then $p\eq1\pmod 6$, $\langle a\rangle_p=(p-1)/6$ and $t=-1/6$. By the above,
\begin{align*}
T_n\l(-\f16\r)\equiv&\frac23-\frac23pq_p(2)-\frac p3H_{\lfloor p/6\rfloor}
\\\eq&\frac23-\frac23pq_p(2)-\frac p3\l(-2q_p(2)-\frac32q_p(3)\r)
\\\eq&\frac23+\frac p2q_p(3)\pmod{p^2}
\end{align*}
and thus
$$\sum_{k=0}^n\frac{\binom{6k}{3k}\binom{3k}k}{(2k+1)432^k}=\frac32T_n\l(-\f16\r)\equiv1+\frac34pq_p(3)=\frac{3^p+1}4\pmod {p^2}.$$

If $a=-1/4$, then $p\eq1\pmod 4$, $\langle a\rangle_p=(p-1)/4$ and $t=-1/4$. By the above,
\begin{align*}
T_n\l(-\f14\r)\equiv&\frac12-pq_p(2)-\frac p2H_{\lfloor p/4\rfloor}
\\\eq&\frac12-pq_p(2)-\frac p2(-3q_p(2))
\\\eq&\frac12+\frac p2q_p(2)\pmod{p^2}
\end{align*}
and thus
$$\sum_{k=0}^n\frac{\binom{4k}{2k}\binom{2k}k}{(2k+1)64^k}=2T_n(-1/4)\equiv1+pq_p(2)=2^{p-1}\pmod {p^2}.$$

{\it Case} II. $\langle a\rangle_p>n$.

If $a=-1/6$, then $p\eq5\pmod 6$,  $\langle a\rangle_p=(5p-1)/6$ and $t=-5/6$. By the above,
\begin{align*}
T_n\l(-\f16\r)\equiv&-\frac23+\frac23pq_p(2)+\frac p3H_{\lfloor {5p}/6\rfloor}
\\\eq&-\frac23+\frac23pq_p(2)+\frac p3\l(-2q_p(2)-\frac32q_p(3)\r)
\\\eq&-\frac23-\frac p2q_p(3)\pmod{p^2}
\end{align*}
and hence $$\sum_{k=0}^n\frac{\binom{6k}{3k}\binom{3k}k}{(2k+1)432^k}=\frac32T_n\l(-\f16\r)\equiv-1-\frac34pq_p(3)=-\frac{3^p+1}4\pmod {p^2}.$$

If $a=-1/4$, then $p\eq3\pmod 4$, $\langle a\rangle_p=(3p-1)/4$ and $t=-3/4$. So
\begin{align*}
T_n\l(-\f14\r)\equiv&-\frac12+pq_p(2)+\frac p2H_{\lfloor{3p}/4\rfloor}
\\\eq&-\frac12+pq_p(2)+\frac p2(-3q_p(2))
\\\eq&-\frac12-\frac p2q_p(2)\pmod{p^2}
\end{align*}
and hence $$\sum_{k=0}^n\frac{\binom{4k}{2k}\binom{2k}k}{(2k+1)64^k}=2T_n\l(-\f14\r)\equiv-1-pq_p(2)=-2^{p-1}\pmod {p^2},$$

The proof of Theorem 1.2 is now complete. \qed

\Ack. The first author would like to thank Dr. Hao Pan for help comments.

\setcounter{conjecture}{0} \end{document}